\documentclass[a4paper, 12pt] {article}
\usepackage[cp1251]{inputenc}
\usepackage[english]{babel}
\usepackage{mathtext}
\usepackage{amsfonts}
\usepackage{amssymb}
\usepackage{amsmath}
\usepackage{xcolor}
\usepackage{cite}
\begin{document}
\thispagestyle{empty}

\noindent{\large\bf ON A REGULARIZATION APPROACH TO THE INVERSE\\
 TRANSMISSION EIGENVALUE PROBLEM}
\\

\noindent {\large\bf S.A.\,Buterin\footnote{Department of
Mathematics, Saratov State University, {\it email:
buterinsa@info.sgu.ru}}, A.E.\,Choque-Rivero\footnote{Instituto de
F\'isica y Matem\'aticas, Universidad Michoacana de San Nicol\'as de
Hidalgo. Edificio C3A, Cd. Universitaria. C. P. 58040 Morelia,
Mich., M\'exico, {\it email: abdon@ifm.umich.mx}} and M.A.
Kuznetsova\footnote{Department of Mathematics, Saratov State
University, {\it email: kuznetsovama@info.sgu.ru}}}
\\

\noindent{\bf Abstract.} We consider the irregular (in the Birkhoff and even the Stone sense) transmission
eigenvalue problem of the form $-y''+q(x)y=\rho^2 y,$ $y(0)=y(1)\cos\rho a-y'(1)\rho^{-1}\sin\rho a=0.$ The
main focus is on the ``most'' irregular case $a=1,$ which is important for applications. The uniqueness
questions of recovering the potential $q(x)$ from transmission eigenvalues were studied comprehensively. Here
we investigate the solvability and stability of this inverse problem. For this purpose, we suggest the
so-called regularization approach, under which there  should first be chosen some regular subclass of
eigenvalue problems under consideration, which actually determines the course of the study and even  the
precise statement of the inverse problem. For definiteness, by assuming $q(x)$ to be a complex-valued
function in $W_2^1[0,1]$ possessing the zero mean value and $q(1)\ne0,$ we study properties of transmission
eigenvalues and prove local solvability and stability of recovering $q(x)$ from the spectrum along with the
value~$q(1).$ In Appendices, we provide some illustrative examples of regular and irregular transmission
eigenvalue problems, and also obtain necessary and sufficient conditions in terms of the characteristic
function for solvability of the inverse problem of recovering an arbitrary real-valued square-integrable
potential $q(x)$ from the spectrum, for any fixed $a\in{\mathbb R}.$


\\

\noindent{\large\bf 1. Introduction}\\

Consider the boundary value problem $R(a,q)$ of the form
\begin{equation}\label{1}
\ell y:=-y''+q(x)y=\lambda y, \quad 0<x<1,
\end{equation}
\begin{equation}\label{2}
y(0)=0, \quad V(y):=y(1)\cos\rho a-y'(1)\frac{\sin\rho a}{\rho}=0,
\end{equation}
where $\rho^2=\lambda$ is the spectral parameter, $q(x)\in L_2(0,1)$ and $a\in{\mathbb R}.$ For $a>0,$ the
problem $R(a,q)$ belongs to the so-called transmission eigenvalue problems. Recently, it has attracted much
attention in connection with the inverse acoustic scattering problem (see \cite{MP, CCM, AGP, CH, CL, WX,
CLM, BYY, BY, XYBY, GP17, Pal, BondB, YB20, WW20} and references therein). A special place among these works
is occupied by studying the {\it inverse} transmission eigenvalue problem, when the potential $q(x)$ is to be
found either on the entire interval $(0,1)$ or on its subinterval from eigenvalues of the problem $R(a,q)$ or
their subset.

The most complete results in the inverse spectral theory are known
for the Sturm--Liouville operator $\ell$ with regular boundary
conditions both in self-adjoint and in non-self-adjoint cases (see,
e.g., \cite{MO, B, Kar, M, HL, L, FY, Mak, B1, BSY, BK}). In
particular, Borg \cite{B} proved that the real-valued potential
$q(x)\in L_2(0,\pi)$ is uniquely determined by specifying the
spectra $\{\lambda_{k,j}\},\;j=0,1,$ of two boundary value problems
${\cal L}_j( q),\;j=0,1,$ for equation (\ref{1}) with one common
boundary condition, for example:
$$
y(0)=y^{(j)}(1)=0,
$$
respectively. For complex-valued potentials, i.e. in the
non-self-adjoint case, this uniqueness result was generalized by
Karaseva \cite{Kar}. It is known that the following asymptotics
hold:
\begin{equation}\label{asympt}
\lambda_{k,j}= \pi^2\Big(k-\frac{j}{2}\Big)^2+\omega+\varkappa_{k,j}, \quad \{\varkappa_{k,j}\}\in l_2, \quad
k\ge1, \quad j=0,1.
\end{equation}
Moreover,
\begin{equation}\label{omega}
\omega=\int_0^1q(x)\,dx.
\end{equation}
Borg \cite{B} also established local solvability and stability of the corresponding inverse problem.
Specifically, the following theorem holds (see also \cite{FY}).

\medskip
{\bf Theorem 1. }{\it For any model real-valued potential $q(x)\in L_2(0,1),$ there exists $\delta>0$ such
that if the arbitrary real sequences $\{\tilde \lambda_{k,j}\}_{k\ge1},\;j=0,1$ satisfy the condition
$$
\Omega:=\sqrt{\sum_{k=1}^\infty\Big(|\lambda_{k,0}-\tilde\lambda_{k,0}|^2+|\lambda_{k,1}-\tilde\lambda_{k,1}|^2\Big)} \le\delta,
$$
then there exists a unique function $\tilde q(x)\in L_2(0,1)$ such that $\{\tilde \lambda_{k,j}\}_{k\ge1}$
are the spectra of the problems ${\cal L}_j(\tilde q),\;j=0,1,$ respectively. Moreover,
$$
\|q-\tilde q\|_2\le C_{q,\delta}\Omega,
$$
where $C_{q,\delta}$ is independent of $\tilde q(x)$ and $\|\cdot\|_\nu:=\|\cdot\|_{L_\nu(0,1)}.$}

\medskip
The original proof of Theorem~1 is also applicable for complex-valued potentials but under the requirement of
simplicity of the spectra. In \cite{BK}, Theorem~1 was generalized for arbitrary multiple spectra: it remains
completely true after replacing all entries of ``real'' with ``complex''. In the self-adjoint case, i.e. when
the function $q(x)$ is real-valued, however, one can prove {\it global} solvability of this inverse problem.
Namely, the following theorem holds (see \cite{MO}).

\medskip
{\bf Theorem 2. }{\it For two arbitrary sequences $\{\lambda_{k,0}\}_{k\ge1}$ and $\{\lambda_{k,1}\}_{k\ge1}$
to be the spectra of the boundary value problems ${\cal L}_0(q)$ and ${\cal L}_1(q),$ respectively, with a
real-valued potential $q(x)\in L_2(0,1),$ it is necessary and sufficient to be real, to have asymptotics
(\ref{asympt}) and to interlace:
\begin{equation}\label{interlace}
\lambda_{k,1}<\lambda_{k,0}<\lambda_{k+1,1}, \quad k\ge1.
\end{equation}}

Unlike the classical Sturm--Liouville problem (when $a=0),$ the boundary conditions (\ref{2}) for $a>0$ can
generally be classified as irregular in the Birkhoff (and even the Stone) sense (see, e.g., \cite{Nai, Fr})
since Green's function of the problem $R(a,q)$ may exponentially grow (see Appendix~A). This results in more
complicated behavior of the spectrum $\{\lambda_k\}$ of $R(a,q).$ Eigenvalues $\lambda_k$ with an account of
multiplicity coincide with zeros of the characteristic function $\Delta(\lambda):=V(S(x,\lambda)),$ where
$y=S(x,\lambda)$ is a solution of equation (\ref{1}) under the initial conditions $S(0,\lambda)=0$ and
$S'(0,\lambda)=1.$ Since $\Delta(\lambda)$ is an entire function of the order not exceeding $1/2,$ according
to Hadamard's factorization theorem, we have
\begin{equation}\label{3}
\Delta(\lambda)=\gamma\Theta(\lambda), \quad \Theta(\lambda)=\lambda^s\prod_{\lambda_k\ne0}\Big(1-\frac{\lambda}{\lambda_k}\Big),
\end{equation}
where $s\ge0$ is the algebraic multiplicity of the zero eigenvalue. Using the transformation operator (see formula (\ref{10-2}) below) for
the solution $S(x,\lambda),$ one can also get the representation
\begin{equation}\label{3a}
\Delta(\lambda)=\frac{\sin\rho(1-a)}{\rho}-\omega\frac{\cos\rho(1-a)}{2\rho^2}+ \int_{a-1}^{a+1} w(t)\frac{\cos\rho t}{\rho^2}\,dt, \quad
w(t)\in L_2(a-1,a+1).
\end{equation}
Unlike the regular case when $a\le0,$ the term $\rho^{-1}\sin\rho(1-a)$ in (\ref{3a}) for $a>0$ is no longer
 the global main part of the
asymptotics for $\Delta(\lambda).$
 Moreover, it even disappears, when $a=1.$

Under the real-valuedness of $q(x),$ however, properties of $R(a,q)$ tend to the self-adjoint case. Namely,
McLaughlin and Polyakov \cite{MP} proved that for real-valued potentials and $a\ne1$ the problem $R(a,q)$ has
infinitely many real eigenvalues $\{\mu_n\}_{n\ge n_0}$ of the form
$$
\mu_n=\frac{\pi^2 n^2}{(1-a)^2}+\frac{\omega}{1-a} +\kappa_n, \quad \{\kappa_n\}\in l_2,
$$
which can always be supplemented by other, possibly nonreal, eigenvalues $\mu_{1},\ldots,\mu_{n_0-1}$ with
the account of multiplicity up to the sequence $\{\mu_n\}_{n\ge1}.$ Moreover, as was illustrated in
\cite{AGP}, the problem $R(a,q)$ may additionally have an infinite number of nonreal eigenvalues. However, if
$a\ne1,$ then the spectrum of $R(a,0),$ obviously, coincides with $\{(1-a)^{-2}\pi^2n^2\}_{n\ge1}.$ In view
of this, the sequence $\{\mu_n\}_{n\ge1}$ was referred to in \cite{BY} as an {\it almost real subspectrum} of
$R(a,q).$ In \cite{MP}, it was proved that specification of $\{\mu_n\}_{n\ge1}$ determines the potential
$q(x)$ uniquely on the subinterval $(0,|a-1|/2),$ if $q(x)$ is known on $(|a-1|/2,1)$   {\it a priori}. In
particular, if $a\ge3,$ then $q(x)$ is determined on the entire interval $(0,1).$ For $a=0$ the corresponding
fact was known as the Hochstadt--Lieberman theorem \cite{HL}. The minimality of the input data
$\{\mu_n\}_{n\ge1}$ for this uniqueness result was established in \cite{BondB}. Moreover, in \cite{BondB}
local solvability and stability of the corresponding inverse problem were proved, having become the first
result dealing with solvability and stability of the inverse transmission eigenvalue problem.

Aktosun and co-authors \cite{AGP} studied the uniqueness of recovering $q(x)$ from the {\it full} spectrum
$\{\lambda_k\}$ of the problem $R(a,q)$ for $a\ge1.$ They reduced the inverse problem to the
 classical inverse Sturm--Liouville problem \cite{B} and proved that $q(x)$ is
uniquely determined by $\{\lambda_k\}$ if $a>1,$ and by $\{\lambda_k\}$ along with the constant $\gamma$ in
(\ref{3}) if $a=1.$ In \cite{BYY}, it was shown that in the case $a=1$ for each nonzero real-valued potential
$q(x)$ one can construct infinitely many different real-valued potentials $\tilde q(x)\in L_2(0,1)$ such that
the corresponding problems $R(1,\tilde q)$ have one and the same spectrum coinciding with the spectrum of
$R(1,q),$ which means the necessity of specifying $\gamma.$  In \cite{BY}, however, it was shown that for
$a>1$ the uniqueness theorem in \cite{AGP} can be improved. Namely, if $a>1,$ then for the unique
determination of the potential it is sufficient to specify only $\{\lambda_k\}\setminus\{\mu_k\}_{k\ge1},$
i.e. the full spectrum with the exception of the entire almost real subspectrum. Moreover,
 even though the authors of \cite{AGP} assumed the real-valuedness of
the potential $q(x),$ their uniqueness results remain true also for complex-valued potentials. In Appendix~B,
we show, in particular, that it holds for $a\le-1$ as well.

The case $a=1$ is exceptional because in general it allows saying almost nothing about the spectrum. For
example, the spectrum of $R(1,0)$ coincides with the entire plane~${\mathbb C},$ while the spectrum of the
problem $R(1,t-1/2)$ is $\{\pi^2 k^2/4+k\varkappa_k\}_{k\ge2},$ where $\{\varkappa_k\}_{k\ge2}$ is a
square-summable sequence. Unlike $R(1,0),$ the problem $R(1,t-1/2)$ obeys some regularization conditions on
the potential, which are stated in the hypothesis of the following theorem.

\medskip
{\bf Theorem 3. }{\it Let $q(x)\in W_2^1[0,1]$ and $q(1)\ne0,$ while $\omega=0.$ Then the spectrum
$\{\lambda_k\}_{k\ge2}$ of the problem $R(1,q)$ has the form
\begin{equation}\label{asympt-1}
\lambda_k=\frac{(\pi k)^2}4+k\varkappa_k, \quad \{\varkappa_k\}\in l_2, \quad k\ge2.
\end{equation}
}

Note that under the hypothesis of Theorem~3 the problem $R(1,q)$ is Stone regular, i.e. its Green's function
polynomially grows as $\lambda\to\infty$ (see Example~A3 in Appendix~A).

In the present paper, we demonstrate the so-called {\it regularization approach}
 that consists of choosing and studying an appropriate regular
subclass of generally speaking irregular eigenvalue problems. The definition of such a class can be given in
terms of some restrictions on the potential $q(x).$ Note that in \cite{XYBY} this idea was used for $a\ne1.$
 Here, however, we apply it to
studying solvability and stability of the inverse problem. For definiteness, we confine ourselves to the
class ${\mathfrak R}$ of problems $R(1,q)$ that is determined by the hypothesis of Theorem~3 and consider the
following inverse problem.

\medskip
{\bf Inverse Problem 1.} Given the spectrum $\{\lambda_k\}_{k\ge2}$ of a problem $R(1,q)\in{\mathfrak R}$
along with the value $\eta:=q(1)/4,$ find the function $q(x).$

\medskip
For our purpose, one can use the complex generalization of Theorem~1 (see \cite{BK}). Therefore, we are able to work in the class of
complex-valued potentials.
  It is more convenient, however, to reduce Inverse Problem~1
     to the problem of
   recovering $q(x)$ from the so-called Cauchy data
(see Section~2). As will be seen below, specification of the value $\gamma$ is equivalent to specification of
$\eta.$ Our main result is the following theorem, which gives
 local solvability and stability of Inverse Problem~1.

\medskip
{\bf Theorem 4. }{\it  Let $\{\lambda_k\}_{k\ge2}$ be the spectrum of
 a certain model problem $R(1,q)$ with a fixed complex-valued
potential $q(x)\in W^1_2[0,1]$ obeying $\eta\ne0$ and
 $\omega=0.$ Then there exists $\delta>0,$ such that for any
  sequence
$\{\tilde\lambda_k\}_{k\ge2}$ and for an arbitrary number
 $\tilde \eta$ satisfying
\begin{equation}\label{11-1}
\Lambda:=|\eta - \tilde \eta| + \sqrt{\sum_{k=2}^\infty
 \frac{|\lambda_k-\tilde\lambda_k|^2}{k^2}} \le\delta,
\end{equation}
there exists a unique problem $R(1,\tilde q)\in{\mathfrak R}$ whose spectrum coincides with the sequence
$\{\tilde\lambda_k\}_{k\ge2}$ and $\tilde q(1)=4\tilde\eta.$ Moreover, the estimate
$$
\|q-\tilde q\|_{W^1_2[0,1]}\le C_{q,\delta}\Lambda
$$
is fulfilled, where $C_{q,\delta}$ is independent of $\tilde q(x),$ and $\|f\|_{W_2^1[a,b]}=\|f\|_{L_2(a,b)}
+\|f'\|_{L_2(a,b)}.$ }

\medskip
Theorem~4, in particular, illustrates the minimality of the input data in Inverse Problem~1. Moreover, it is
the first local solvability and stability result in the inverse transmission eigenvalue problem for
complex-valued potentials. The proof of Theorem~4 is constructive.

We note that the suggested regularization approach is vital for finding conditions for solvability of an
inverse problem in terms of the spectrum. However, sometimes one can alternatively formulate conditions for
solvability in terms of the characteristic function. For example, in the recent work \cite{WW20} this was
suggested for the problem of recovering an arbitrary real-valued potential $q(x)\in L_2(0,1)$ from the
spectrum of $R(1,q).$ But, unfortunately, the corresponding theorem (Theorem 4.1) contains a mistake. In
Appendix~B, by using the results of \cite{BYY} we correct the mentioned mistake in \cite{WW20} and extend
this result to all other real values of the parameter~$a.$

The paper is organized as follows. In the next section, we provide some auxiliary results and obtain an
algorithm for solving Inverse Problem~1. The proof of Theorem~4 is given in Section~3. In Appendix~A, we
provide several examples of regular and irregular problems $R(a,q).$ In Appendix~B, we obtain necessary and
sufficient conditions in terms of the characteristic function for solvability of the inverse problem of
recovering an arbitrary real-valued potential $q(x)\in L_2(0,1)$ from the spectrum of the problem $R(a,q)$
for any real $a.$

Throughout the paper, one and the same symbol $C_{q,\delta}$ denotes {\it different} positive constants in estimates, which depend only on
$q(x)$ and $\delta.$
\\

\noindent{\large\bf 2. Constructive solution of the inverse problem}\\

We start with the following well-known representation (see, e.g., \cite{M}):
\begin{equation}\label{10-2}
S(x, \lambda) = \frac{\sin \rho x}{\rho} + \int_0^x
K(x,t)\frac{\sin\rho t}\rho\,dt, \quad 0\le x\le1,
\end{equation}
where $K(x,t)$ is a continuous function and $K(x,0)=0.$ More precisely,
 by virtue
 of formulae (1.2.9) and
(1.2.17) in \cite{M}, after the odd continuation
\begin{equation}\label{odd}
K(x,t):=-K(x,-t), \quad -1\le-x\le t<0,
\end{equation}
and then the continuation by zero outside the triangle $0\le|t|\le
x\le1,$ the kernel $K(x,t)$ will satisfy the integral equation
\begin{equation}\label{K0}
K(x,t)=\frac12\int_{\frac{x-t}{2}}^{\frac{x+t}{2}}q(\tau)\,d\tau
+\frac12\int_0^xq(\tau)\,d\tau\int_{t-(x-\tau)}^{t+(x-\tau)} K(\tau,\xi)\,d\xi, \quad 0\le|t|\le x\le1.
\end{equation}
Note that the domain of integration in the double integral in (\ref{K0}) includes subdomains on which
$|\xi|>\tau$ and, hence, $K(\tau,\xi)$ may possess first-order discontinuities. In order to remove them, it
is sufficient to rewrite the equation (\ref{K0}) in the following
 equivalent form
 \begin{align}
K(x,t)=&\frac12\int_{\frac{x-t}{2}}^{\frac{x+t}{2}}q(\tau)\,d\tau
+\frac12\int_\frac{x-t}2^xq(\tau)\,d\tau\int_{t-(x-\tau)}^\tau
K(\tau,\xi)\,d\xi \nonumber\\
&
 -\frac12\int_\frac{x+t}2^xq(\tau)\,d\tau\int_{t+x-\tau}^\tau K(\tau,\xi)\,d\xi,
\quad 0\le|t|\le x\le1, \label{K}
\end{align}
in the right-hand side of which the inequality $|\xi|\le\tau$ automatically holds.

Alternatively, by substituting~(\ref{10-2}) directly into equation (\ref{1}) and integrating by parts, one
can show that after the continuation (\ref{odd}) the kernel $K(x,t)$ becomes a solution of the following
Goursat problem (see \cite{L}):
\begin{equation}\label{Partial}
K_{xx}(x,t)-K_{tt}(x,t)=q(x)K(x,t), \quad 0<|t|<x\le1,
\end{equation}
\begin{equation}\label{Goursat}
K(x,\pm x)=\pm\frac12\int_0^x q(t)\,dt.
\end{equation}
For potentials $q(x)\notin W_1^1[0,1],$ we emphasize that the second partial derivatives in (\ref{Partial})
do not exist in the usual sense. Therefore, finding and studying the kernel $K(x,t)$ are more convenient
directly via the integral equation (\ref{K}), which is equivalent to the Goursat problem (\ref{Partial}),
(\ref{Goursat}) and can be derived independently (see \cite{M}).

Let $j\in\{0,1\}.$ Eigenvalues of the problem ${\cal L}_j(q)$ coincide
 with zeros if its characteristic function
$\Delta_j(\lambda):=S^{(j)}(1,\lambda).$ Integrating by parts and differentiating
 in (\ref{10-2}), we obtain
\begin{equation}\label{1-3}
\Delta_0(\lambda)=\frac{\sin\rho}{\rho} -\omega\frac{\cos\rho}{2\rho^2}+ \int_0^1 w_0(t)\frac{\cos\rho t}{\rho^2}\,dt, \;\;
\Delta_1(\lambda)=\cos\rho +\omega\frac{\sin\rho}{2\rho}+ \int_0^1 w_1(t)\frac{\sin\rho t}{\rho}\,dt,
\end{equation}
where $w_0(t) = K_t(1, t)$ and $w_1(t) = K_x(x, t)|_{x=1}.$
  Thus, the kernel $K(x,t)$
 is also a solution of the Cauchy problem for the
equation (\ref{Partial}) along with the initial conditions
\begin{equation}\label{Cauchy}
K(1,t)=\int_0^t w_0(\tau)\,d\tau, \quad K_x(x,t)|_{x=1}=w_1(t), \quad t\in[-1,1],
\end{equation}
where $w_j(t)=(-1)^jw_j(-t)$ for $t\in[-1,0)$ and $j=0,1.$ Within this context, the ordered pair of functions
$\{w_0,w_1\}$ is sometimes referred to as {\it Cauchy data} related to the potential $q(x).$ Thus, after
assuming $\omega$ to be fixed, Borg's statement of the inverse
 problem, which consists of recovering $q(x)$ from the spectra
 $\{\lambda_{k,0}\}$ and $\{\lambda_{k,1}\},$ is equivalent to the following
inverse problem from the Cauchy data.

\medskip
{\bf Inverse Problem 2.} Given the functions $w_0(x)$ and $w_1(x),$ find the function $q(x)$ such that the solution $K(x,t)$ of the Goursat
problem (\ref{Partial}), (\ref{Goursat}) satisfies the conditions (\ref{Cauchy}).

\medskip
 By substituting (\ref{1-3}) into
  $\Delta(\lambda)=V(S(x,\lambda))=\Delta_0(\lambda)\cos\rho a-\rho^{-1}\Delta_1(\lambda)\sin\rho a,$
  we arrive
at the representation (\ref{3a}), in which we have
\begin{equation} \label{wa}
w(t) = \frac{1}{2} \left\{\begin{array}{cc}
w_0(a-t)-w_1(a-t), & t\in[a-1,a],\\[2mm]
w_0(t-a)+w_1(t-a), & t\in(a,a+1].
\end{array} \right.
\end{equation}
For $a=1,$ it takes the form
\begin{equation} \label{Delta}
\Delta(\lambda) = -\frac\omega{2\rho^2} +\int_0^2 w(t) \frac{\cos \rho t}{\rho^2} \, dt,
\end{equation}
where
\begin{equation} \label{w}
w(t) = \frac{1}{2}
\left\{\begin{array}{cc}
w_0(1-t) - w_1(1-t) , & t \in [0, 1],\\[2mm]
w_0(t-1) + w_1(t - 1), & t \in (1, 2].
\end{array} \right.
\end{equation}
In what follows, we assume $R(1,q)\in{\mathfrak R},$ i.e. the function $q(x)$
 obeys the hypothesis of Theorem~3. Then, by using (\ref{K}), one can
show that $w(x)\in W^1_2[0,2].$ Indeed, since $q(x)\in W_2^1[0,1],$
we have $w_0(x),\,w_1(x)\in W_2^1[0,1].$ Thus, it remains to note
that $w(1-0) = w(1+0)$ since $w_1(0)=0.$ Let us calculate $w(2).$
 By differentiating (\ref{K}), we arrive at
 \begin{align*}
 K_x(x,t)=&\frac14\Big(q\Big(\frac{x+t}2\Big)
-q\Big(\frac{x-t}2\Big)\Big) +\frac12\int_\frac{x-t}2^xq(\tau)
K(\tau,t+\tau-x)\,d\tau\\
&+\frac12\int_\frac{x+t}2^xq(\tau) K(\tau,t+x-\tau)\,d\tau,
 \end{align*}
 \begin{align}
K_2(x,t):=K_t(x,t)&=\frac14\Big(q\Big(\frac{x+t}2\Big) +q\Big(\frac{x-t}2\Big)\Big)
-\frac12\int_\frac{x-t}2^xq(\tau)
K(\tau,t+\tau-x)\,d\tau \nonumber\\
 &
+\frac12\int_\frac{x+t}2^xq(\tau) K(\tau,t+x-\tau)\,d\tau.
 \label{Kt}
\end{align}
Hence, in particular,
$$
w_0(1)=\frac{q(1)+q(0)}4 -\frac12\int_0^1
q(\tau)K(\tau,\tau)\,d\tau, \quad w_1(1)=\frac{q(1)-q(0)}4
+\frac12\int_0^1 q(\tau)K(\tau,\tau)\,d\tau.
$$
By virtue of (\ref{w}), we finally get
\begin{equation}\label{2.1}
w(2)=\frac{w_0(1)+w_1(1)}2= \frac{q(1)}4 = \eta.
\end{equation}
Thus, we have
\begin{equation} \label{rec_w}
w(x) = \eta + \int_x^2 v(t) \, dt,  \quad v(t) = -w'(t).
\end{equation}
By integrating by parts in (\ref{Delta}) and by taking into account that $\omega=0$ for $R(1,q)\in{\mathfrak
R},$ we get
\begin{equation} \label{partsDelta}
\Delta(\lambda) =  \eta \frac{\sin 2 \rho}{\rho^3} + \int_0^2 v(t)  \frac{\sin \rho t}{\rho^3} \, dt, \quad \eta\ne0, \quad v(t) \in L_2(0,
2).
\end{equation}
By the standard approach involving Rouch\'e's theorem
 (see, e.g., \cite{FY}), one can show that any entire function $\Delta(\lambda)$ of
the form (\ref{partsDelta}) has infinitely many zeros $\lambda_k,\,k\ge2,$ of the form (\ref{asympt-1}),
which gives the assertion of Theorem~3. Moreover, using Hadamard's factorization theorem, by the standard
approach (see, e.g., \cite{FY}) one can prove that the function $\Delta(\lambda)$ is determined by its zeros
along with the constant $\eta$ uniquely. Moreover, the following formula holds:
\begin{equation}\label{char2}
\Delta(\lambda)=-\frac{8\eta}{\pi^2}\prod_{k=2}^\infty \frac{4(\lambda_k-\lambda)}{(\pi k)^2}.
\end{equation}
Conversely, the following lemma can be obtained as a corollary from
  Lemma~3.3 in \cite{But07}.

\medskip
{\bf Lemma 1. }{\it For any complex sequence $\{\lambda_k\}_{k\ge2}$
 of the form (\ref{asympt-1}), the function $\Delta(\lambda)$
  determined
by formula (\ref{char2}) with some $\eta\ne0$ has the form
(\ref{partsDelta}) with some function $v(t)\in L_2(0,2).$}

\medskip
The next lemma, being a corollary from Lemma~1 in \cite{But20}, gives uniform stability of recovering the function $v(t)$ from given zeros
$\{\lambda_k\}_{k\ge2}$ of the function $\Delta(\lambda)$ along with the value~$\eta.$

\medskip
{\bf Lemma 2. }{\it For any $r>0,$ there exists $C_r>0$ such that
$$
\|v-\tilde v\|_{L_2(0,2)} \le C_r\Lambda
$$
as soon as $|\eta|\le r$ (or alternatively, $|\tilde\eta|\le r\!\!$)
and
$$
\sum_{k=2}^\infty \frac{|4\lambda_k-(\pi k)^2|^2}{k^2} \le r, \quad \sum_{k=2}^\infty
\frac{|4\tilde\lambda_k-(\pi k)^2|^2}{k^2} \le r.
$$
Here $\Lambda$ is determined in (\ref{11-1}), while the function $\tilde v(x)$ is
 determined by the relation
\begin{equation}\label{char_ti}
\tilde\Delta(\lambda):=-\frac{8\tilde\eta}{\pi^2}\prod_{k=2}^\infty \frac{4(\tilde\lambda_k-\lambda)}{(\pi
k)^2} = \tilde \eta \frac{\sin 2 \rho}{\rho^3} + \int_0^2 \tilde v(x)  \frac{\sin \rho x}{\rho^3} \, dx.
\end{equation}
}

Now we are in the position to provide an algorithm for solving
Inverse Problem~1. Fix a model problem $R(1,q)\in{\mathfrak R}$ with
the spectrum $\{\lambda_k\}_{k\ge2}.$ Let an arbitrary nonzero
complex number~$\tilde\eta$ and a complex sequence
$\{\tilde\lambda_k\}_{k\ge2}$ be given that obey inequality
(\ref{11-1}) with a sufficiently small $\delta>0.$ Thus, the
corresponding potential $\tilde q(x)$ can be found by the following
algorithm.

\medskip
{\bf Algorithm 1. }{\it (i) Construct the function $\tilde v(x)$ by the formula
\begin{equation}\label{v-ti}
\tilde v(x)=\frac{\pi^3}8\sum_{k=1}^\infty k^3\tilde\Delta\Big(\frac{\pi^2k^2}4\Big)\sin\frac{\pi kx}2,
\end{equation}
where the function $\tilde\Delta(\lambda)$ is determined by the first
 equality in (\ref{char_ti}).

(ii) Calculate the functions $\tilde w_0(x)$ and $\tilde w_1(x)$ by the formulae
\begin{equation}\label{w-j_ti}
\tilde w_j(x)=\tilde w(1+x)+(-1)^j\tilde w(1-x), \quad x\in[0,1], \quad j=0,1,
\end{equation}
where the function $\tilde w(x)$ is determined by the formula
\begin{equation}\label{w_ti}
\tilde w(x) = \tilde\eta + \int_x^2 \tilde v(t)\,dt.
\end{equation}

(iii) For $j=0,1,$ find zeros $\{\tilde\lambda_{k,j}\}_{k\ge1}$ of the function $\tilde\Delta_j(\lambda),$ where
\begin{equation}\label{1-3-ti}
\tilde\Delta_0(\lambda)=\frac{\sin\rho}{\rho} + \int_0^1 \tilde w_0(x)\frac{\cos\rho x}{\rho^2}\,dx, \quad
\tilde\Delta_1(\lambda)=\cos\rho + \int_0^1 \tilde w_1(x)\frac{\sin\rho x}{\rho}\,dx.
\end{equation}

(iv) Put $\tilde q(x)=q(x)+r(x),$ where $r(x)$ is a solution of the Borg equation~(38) in \cite{BK}.}

\medskip
Using Lemma~2 along with a $W_2^1$-analogue of Lemma~4.6 in \cite{YBX} as well as Theorem~1 for
complex-valued potentials (see Theorem~3 in \cite{BK}), one can show that for sufficiently small $\delta>0$
the Borg equation in step~(iv) of Algorithm~1 is uniquely solvable. At the same time, the following example
shows, in particular, that the choice of sufficiently small $\delta>0$ is important.

\medskip
{\bf Example 1.} Let $\tilde\lambda_k=\pi^2k^2/4,\,k\ge2,$ then formulae (\ref{char_ti}) and (\ref{v-ti})
give
$$
\tilde\Delta(\lambda)=-\frac{8\tilde\eta}{\pi^2}\prod_{k=2}^\infty \Big(1-\frac{4\lambda}{(\pi k)^2}\Big) =
\tilde \eta \frac{\sin 2 \rho}{\rho^3} + \int_0^2 \tilde v(x)  \frac{\sin \rho x}{\rho^3} \, dx, \quad \tilde
v(x)=-\frac{\pi\tilde\eta}2\sin\frac{\pi x}2.
$$
Then, by using formulae (\ref{w-j_ti}) and (\ref{w_ti}) we calculate
$$
\tilde w(x)=-\tilde\eta\cos\frac{\pi x}2, \quad \tilde w_j(x)=-\tilde\eta\Big(\cos\frac\pi2(1+x)
+(-1)^j\cos\frac\pi2(1-x)\Big) =2j\tilde\eta\sin\frac{\pi x}2,
$$
where $j=0,1,$ which along with (\ref{1-3-ti}) give $\tilde\Delta_0(\lambda)=\rho^{-1}\sin\rho$ and
$$
\tilde\Delta_1(\lambda)=\cos\rho +\frac{\tilde\eta}\rho\int_0^1\Big(\cos\Big(\frac\pi2-\rho\Big) x
-\cos\Big(\frac\pi2+\rho\Big) x\Big)dx \qquad\qquad\qquad\qquad\qquad\qquad
$$
$$
=\cos\rho +\frac{\tilde\eta}\rho\Big(\Big(\frac\pi2-\rho\Big)^{-1}\sin\Big(\frac\pi2-\rho\Big)
-\Big(\frac\pi2+\rho\Big)^{-1}\sin\Big(\frac\pi2+\rho\Big)\Big) \qquad\qquad
$$
$$
\quad\qquad =\cos\rho +\tilde\eta\Big(\frac{\pi^2}4-\lambda\Big)^{-1}\Big(\sin\Big(\frac\pi2-\rho\Big)
+\sin\Big(\frac\pi2+\rho\Big)\Big) =(\pi^2-4\lambda+8\tilde\eta)\frac{\cos\rho}{\pi^2-4\lambda}.
$$
Thus, the third step of Algorithm~1 gives
$$
\tilde\lambda_{k,0}=\pi^2k^2, \;\; k\ge1, \quad \tilde\lambda_{1,1}=\frac{\pi^2}4+2\tilde\eta, \;\;
\tilde\lambda_{k,1}=\pi^2\Big(k-\frac12\Big)^2, \;\; k\ge2.
$$
According to Theorem~2, there exists a real-valued potential $\tilde q(x)$ such that the constructed
sequences $\{\tilde\lambda_{k,0}\}_{k\ge1}$ and $\{\tilde\lambda_{k,1}\}_{k\ge1}$ are the spectra of the
problems ${\cal L}_0(\tilde q)$ and ${\cal L}_1(\tilde q),$ respectively, if and only if
$\tilde\eta<3\pi^2/8.$ Thus, taking (\ref{2.1}) into account, one can see that, solvability of Inverse
Problem~1 with the input data, consisting of the sequence $\{\pi^2k^2/4\}_{k\ge2}$ along with the number
$\tilde\eta,$ in the class of real-valued potentials is equivalent to
$\tilde\eta\in(-\infty,0)\cup(0,3\pi^2/8).$

Taking, for example, $\eta=\pi^2/4$ and $\lambda_k=\pi^2k^2/4=\tilde\lambda_k,$ $k\ge2,$ as the model input
data, one can see that $\delta$ in (\ref{11-1}) should be less than $\pi^2/8.$ Otherwise, it would admit the
value $\Lambda=\pi^2/8$ allowing $\tilde\eta$ to be equal to $3\pi^2/8,$ i.e.
$\tilde\lambda_{1,1}=\pi^2=\tilde\lambda_{1,0},$ which leads to nonexistence of $\tilde q(x),$ since the
problems ${\cal L}_0(\tilde q)$ and ${\cal L}_1(\tilde q)$ cannot possess common eigenvalues.

\medskip
For proving Theorem~4, it is convenient, however, to replace steps (iii) and (iv) in Algorithm~1 with direct
recovering the potential $\tilde q(x)$ from the Cauchy data $\{\tilde w_0,\tilde w_1\}.$ Recently, Bondarenko
\cite{Bond20} proved the following theorem, which gave local solvability and stability of Inverse Problem~2
for complex-valued potentials (see Theorem~5.1 in Appendix of \cite{Bond20}).

\medskip
{\bf Theorem 5. }{\it For each complex-valued potential $q(x)\in L_2(0,1),$ there exists $\varepsilon>0$ such
that for any functions $\tilde w_j(x)\in L_2(0,1),\,j=0,1,$ satisfying the estimate
\begin{equation} \label{CauDel}
\Xi:=\max_{j=0,1}\|w_j-\tilde w_j\|_2\le\varepsilon,
\end{equation}
there exists a unique function $\tilde q(x)\in L_2(0,1)$ such that $\int_0^1q(x)\,dx=\int_0^1\tilde q(x)\,dx$ and $\{\tilde w_0,\tilde
w_1\}$ are the Cauchy data for $\tilde q(x).$ Moreover, the following estimate holds:
\begin{equation} \label{CauStab}
\|q-\tilde q\|_2\le C_{q,\varepsilon}\Xi.
\end{equation}
Here the pair $\{w_0,w_1\}$ is the Cauchy data related to the potential $q(x).$}

\medskip
In the next section, leaning on this result we give the proof of Theorem~4. One of the main technical
difficulties is connected with our dealing with $W_2^1$-potentials.
  It will be seen that for our purpose, however, there is no need to
   derive any $W_2^1$-analogue of Theorem~5.
\\

\noindent{\large\bf 3. Proof of Theorem~4}\\

Fix a problem $R(1,q)\in{\mathfrak R}$ with the spectrum $\{\lambda_k\}_{k\ge2}.$ Let us be given with a certain nonzero complex number
$\tilde\eta$ and some complex sequence $\{\tilde\lambda_k\}_{k\ge2}$ for which the value $\Lambda$
 determined in (\ref{11-1}) is finite.
 It is then easy to see that $\tilde\lambda_k=(\pi k)^2/4+k\tilde\varkappa_k,$ where $\{\tilde\varkappa_k\}\in l_2.$
  By virtue of Lemma~1,
there exists a unique function $\tilde v(x)$ for which representation (\ref{char_ti}) is fulfilled. Determine
the functions $\tilde w_0(x)$ and $\tilde w_1(x)$ by formula (\ref{w-j_ti}), where the function $\tilde w(x)$
is determined by formula (\ref{w_ti}). By using (\ref{w}) combined with (\ref{w-j_ti}) and (\ref{rec_w})
along with (\ref{w_ti}), it is easy to estimate $\|w_j-\tilde w_j\|_{W_2^1[0,1]} \le\sqrt2\|w-\tilde
w\|_{W_2^1[0,2]}$ for $j=0,1$ and $\|w-\tilde w\|_{W_2^1[0,2]} \le \sqrt2|\eta-\tilde\eta|
+(\sqrt2+1)\|v-\tilde v\|_{L_2(0,2)},$ respectively. By combining these estimates, we get
$$
\|w_j-\tilde w_j\|_{W_2^1[0,1]} \le 2|\eta-\tilde\eta| +(2+\sqrt2)\|v-\tilde v\|_{L_2(0,2)}, \quad j=0,1,
$$
which along with Lemma~2 imply the estimate
\begin{equation}\label{w_hat}
\|w_j-\tilde w_j\|_{W_2^1[0,1]} \le C_{q,\delta}\Lambda, \quad j=0,1,
\end{equation}
as soon as inequality~(\ref{11-1}) is fulfilled. Thus, by virtue of Theorem~5, for sufficiently small
$\delta>0,$ inequality~(\ref{11-1}) implies the existence of a unique potential $\tilde q(x)\in L_2(0,1)$
with the Cauchy data $\{\tilde w_0,\tilde w_1\}.$ Moreover, by virtue of (\ref{CauDel})--(\ref{w_hat}), we
have the estimate
\begin{equation}\label{3.0}
\|q-\tilde q\|_2\le C_{q,\delta}\Lambda.
\end{equation}

Furthermore, since $\tilde w_0(x),\tilde w_1(x)\in W_2^1[0,1],$ we
have $\tilde q(x)\in W_2^1[0,1].$ Indeed, this can be easily
obtained as a consequence from the corollary to Theorem~1.5.1 in
\cite{M}. It is easy to see that the corresponding problem
$R(1,\tilde q)$ belongs to the class ${\mathfrak R}$ and has the
spectrum $\{\tilde\lambda_k\}_{k\ge2}.$
 Moreover, by virtue of
(\ref{2.1}), we have $ \tilde q(1) = 4 \tilde\eta.$

Thus, for finishing the proof of Theorem~4 it remains to establish the estimate
\begin{equation}\label{q'stab}
\|q'-\tilde q'\|_2\le C_{q,\delta}\Lambda.
\end{equation}

 We agree that if some symbol $\alpha$ denotes an object related to
 the potential $q(x),$ then this symbol with tilde $\tilde\alpha$
 denotes
 the analogous object corresponding to $\tilde q(x),$ and
  $\hat\alpha:=\alpha-\tilde\alpha.$ The subsequent arguments
   partially repeat those
 in Borg's method (see \cite{BK}).

Since $\ell S(x,\lambda)=\lambda S(x,\lambda)$ and $\tilde\ell\tilde S(x,\lambda)=\lambda\tilde S(x,\lambda),$ we get
\begin{equation}\label{3.1}
\int_0^1 \hat q(x) S(x,\lambda)\tilde S(x,\lambda)\,dx =\tilde S(1,\lambda)S'(1,\lambda) -\tilde S'(1,\lambda)S(1,\lambda)
=\tilde\Delta_0(\lambda)\Delta_1(\lambda) -\tilde\Delta_1(\lambda)\Delta_0(\lambda).
\end{equation}
Put
\begin{equation}\label{3.2}
\varphi(x,\lambda):=1-2\lambda S(x,\lambda)\tilde S(x,\lambda) =\cos2\rho x +\int_0^x Q(x,t)\cos2\rho t\,dt,
\end{equation}
where $Q(x,t)$ is a continuous function. Moreover, by substituting (\ref{10-2}) into (\ref{3.2}) and
using~(\ref{odd}), one can
 calculate
\begin{equation}\label{3.2.1}
Q(x,t)=2\Big(K(x,2t-x) +\tilde K(x,2t-x) +\int_{2t-x}^x K(x,\tau)\tilde K(x,2t-\tau)\,d\tau\Big),\;\; 0\le t\le x\le1.
\end{equation}
 By substituting (\ref{3.2}) into (\ref{3.1}) and taking into account the zero mean value of $\hat q(x),$ we get
\begin{equation}\label{3.3}
\int_0^1 \hat q(x) \varphi(x,\lambda)\,dx =2\lambda\Big(\Delta_0(\lambda)\tilde\Delta_1(\lambda)
-\tilde\Delta_0(\lambda)\Delta_1(\lambda)\Big).
\end{equation}
It is easy to show that the function $Q(x,t)$ is continuous and possesses square-integrable partial
derivatives $Q_x(x,t)$ and $Q_t(x,t)$ on the triangle $0<t<x<\pi,$ if $q(x),\tilde q(x)\in L_2(0,1),$ while
under our standing condition $q(x),\tilde q(x)\in W_2^1[0,1],$ the kernel $Q(x,t)$ acquires an additional
degree of smoothness. By
 integrating in (\ref{3.3}) by parts and multiplying with $2\rho,$ we arrive at
\begin{equation}\label{3.4}
\int_0^1 \hat q'(x) \phi(x,\lambda)\,dx =\omega_1\Big(\frac{2\rho}\pi\Big)-\omega_2\Big(\frac{2\rho}\pi\Big),
\end{equation}
where $\omega_1(2\rho/\pi)=\hat q(1)\phi(1,\lambda)$ and $\omega_2(2\rho/\pi)=4\rho^3( \Delta_0(\lambda)\tilde\Delta_1(\lambda)
-\tilde\Delta_0(\lambda)\Delta_1(\lambda)),$ while
\begin{equation}\label{3.5}
\phi(x,\lambda)=2\rho\int_0^x\varphi(t,\lambda)\,dt =\sin2\rho x +\int_0^x U(x,t)\sin2\rho t\,dt,
\end{equation}
\begin{equation}\label{3.6}
U(x,t)=-\frac{d}{dt}\int_t^x Q(\tau,t)\,d\tau =Q(t,t)-\int_t^x Q_t(\tau,t)\,d\tau.
\end{equation}
Thus, the functional sequence $\{\phi(x,(\pi n)^2/4)\}_{n\in{\mathbb N}}$ is a Riesz basis in $L_2(0,1).$
Hence, formula (\ref{3.4}) implies the estimate
\begin{equation}\label{3.7}
\|\hat q'\|_2 \le A\Big(\|\{\omega_1(n)\}_{n\in{\mathbb N}}\|_{l_2} +\|\{\omega_2(n)\}_{n\in{\mathbb N}}\|_{l_2}\Big),
\end{equation}
where, according to (\ref{3.5}), we have
$$
A=\sqrt2\|(I+U^*)^{-1}\|=\sqrt2\|(I+U)^{-1}\|, \quad Uf=\int_0^x U(x,t)f(t)\,dt,
$$
while $I$ is the identity operator and $\|\cdot\|:=\|\cdot\|_{L_2(0,1)\to L_2(0,1)}$ (see, e.g., Section 1.8.5 in \cite{FY}).
 Furthermore, by
virtue of Lemma~1 in \cite{ButMal18}, we have the estimate
$$
\|(I+U)^{-1}\|\le1+\|(I+U)^{-1}-I\|\le1+ F(\|U(\,\cdot\,,\,\cdot\,)\|_{L_2((0,1)^2)}), \quad
F(x)=x+\sum_{k=0}^\infty\frac{x^{k+2}}{\sqrt{k!}}.
$$
On the other hand, by solving the integral equation (\ref{K}) with the method of successive approximations
(see, e.g., Theorem~1.2.2 in \cite{M}), one can get the estimate
$$
|K(x,t)|\le \|q\|_1\exp(\|q\|_1), \quad 0\le t\le x\le1,
$$
which along with (\ref{Kt}) yield
$$
\|K_2(\,\cdot\,,\,\cdot\,)\|_{L_2((0,1)^2)}\le \|q\|_2 +\|q\|_1^2\exp(\|q\|_1).
$$
Thus, by using (\ref{11-1}), (\ref{odd}), (\ref{3.0}), (\ref{3.2.1})
and (\ref{3.6}), we get
$$
\|U(\,\cdot\,,\,\cdot\,)\|_{L_2((0,1)^2)}\le C_{q,\delta}, \quad \|U(1,\,\cdot\,)\|_2\le C_{q,\delta}.
$$
Hence, in (\ref{3.7}) we have
\begin{equation}\label{3.7.1}
A\le C_{q,\delta},
\end{equation}
and it remains to prove the estimates
\begin{equation}\label{3.8}
\|\{\omega_j(n)\}_{n\in{\mathbb N}}\|_{l_2}\le C_{q,\delta}\Lambda, \quad j=1,2.
\end{equation}
For $j=1,$ we get
$$
\|\{\omega_1(n)\}_{n\in{\mathbb N}}\|_{l_2}\le |\hat q(1)|\cdot\|\{\phi(1,(\pi n)^2/4)\}_{n\in{\mathbb N}}\|_{l_2}
=4\sqrt2\|U(1,\,\cdot\,)\|_2|\hat\eta| \le C_{q,\delta}\Lambda.
$$
For $j=2,$ we have
$$
\frac1{4\rho^3}\omega_2\Big(\frac{2\rho}\pi\Big) =\Delta_0(\lambda)\tilde\Delta_1(\lambda) -\tilde\Delta_0(\lambda)\Delta_1(\lambda)
=\hat\Delta_0(\lambda)\Delta_1(\lambda) -\Delta_0(\lambda)\hat\Delta_1(\lambda)
$$
where, by using (\ref{1-3}) with $\omega=0$ and (\ref{1-3-ti}) we
obtain
 \begin{align*}
\Delta_0(\lambda)=&\frac{\sin\rho}\rho+\int_0^1 w_0(t)\frac{\cos\rho
t}{\rho^2}\,dt, \quad \Delta_1(\lambda)=\cos\rho+\int_0^1
w_1(t)\frac{\sin\rho t}{\rho}\,dt,
\\
\hat\Delta_0(\lambda)=&\int_0^1\hat w_0(t)\frac{\cos\rho
t}{\rho^2}\,dt =\hat w_0(1)\frac{\sin\rho}{\rho^3}-\int_0^1\hat
w_0'(t)\frac{\sin\rho t}{\rho^3}\,dt,
 \\
\hat\Delta_1(\lambda)=&\int_0^1\hat w_1(t)\frac{\sin\rho
t}{\rho}\,dt =-\hat w_1(1)\frac{\cos\rho}{\rho^2} +\int_0^1\hat
w_1'(t)\frac{\cos\rho t}{\rho^2}\,dt.
 \end{align*}
 Therefore, we have
 \begin{align*}
\omega_2(n)=&-4\int_0^1\hat w_0'(t)\sin\frac{\pi nt}2\,dt
\Big((-1)^\frac{n}2 +\frac2{\pi n}\int_0^1 w_1(t)\sin\frac{\pi
nt}2\, dt\Big)\\
 &-4\int_0^1 w_0(t)\cos\frac{\pi nt}
2\,dt \int_0^1\hat w_1(t)\sin\frac{\pi nt}2\,dt
 \end{align*}
for even $n,$ and
 \begin{align*}
\omega_2(n)=&4\int_0^1\hat w_0(t)\cos\frac{\pi nt}2\,dt
\int_0^1w_1(t)\sin\frac{\pi nt}2\,dt\\
 &-4\Big((-1)^\frac{n-1}2 +\frac2{\pi n}\int_0^1 w_0(t)\cos\frac{\pi nt}2\,dt\Big) \int_0^1\hat w_1'(t)\cos\frac{\pi
nt}2\,dt
 \end{align*}
for odd $n.$ Hence, we get the estimates
$$
\|\{\omega_2(2n)\}_{n\in{\mathbb N}}\|_{l_2}\le C_q(\|\hat w_0'\|_2
+\|\hat w_1\|_2), \quad \|\{\omega_2(2n-1)\}_{n\in{\mathbb
N}}\|_{l_2}\le C_q(\|\hat w_0\|_2 +\|\hat w_1'\|_2).
$$
Thus, by virtue of (\ref{w_hat}), we arrive at the estimate (\ref{3.8}) also for $j=2.$ According to
(\ref{3.7})--(\ref{3.8}), we have (\ref{q'stab}), which finishes the proof. $\hfill\Box$
\\

\noindent{\large\bf Appendix A}
\\

Here we provide several illustrative examples of both regular and irregular problems $R(a,q).$ Denote by
$G(x,t,\lambda)$ the Green's function of $R(a,q),$ which is determined by the formula
$$
y(x)=\int_0^1G(x,t,\lambda)f(t)\,dt,
$$
where $y(x)$ is the solution of the boundary value problem
$$
\ell y=\lambda y+f(x), \quad 0<x<1, \quad y(0)=V(y)=0,
 \quad f(x)\in L_2(0,1).
$$
In accordance with the classical direct spectral theory of ordinary differential operators (see, e.g.,
\cite{Nai,Fr}), we refer to the problem $R(a,q)$ as {\it Birkhoff regular}, if it possesses a Green's
function and there exist expanding contours $\{\lambda:|\lambda|=r_k\},$ where $r_k\to\infty$ as
$k\to\infty,$ on which the estimate
\begin{equation}\label{a.1}
G(x,t,\lambda)=O(\lambda^\theta), \quad \lambda\to\infty,
\end{equation}
is fulfilled for $\theta=-1/2.$ If the Green's function exists and estimate (\ref{a.1}) holds for at least
some finite $\theta,$ then the problem $R(a,q)$ is referred to as {\it Stone regular}. By substitution, it is
easy to check that the function $G(x,t,\lambda)$ has the form
\begin{equation}\label{a.1.1}
G(x,t,\lambda)=-\frac1{\Delta(\lambda)}\left\{
\begin{array}{l} \psi(x,\lambda)S(t,\lambda), \quad t\le x,\\[3mm]
S(x,\lambda)\psi(t,\lambda), \quad t\ge x,
\end{array}\right.
\end{equation}
where $\psi(x,\lambda)$ is a solution of equation (\ref{1}) under the initial conditions
$$
\psi(1,\lambda)=\frac{\sin\rho a}\rho, \quad \psi'(1,\lambda)=\cos\rho a.
$$
The following asymptotics holds (see \cite{BondB}):
\begin{equation}\label{a.2}
\psi(x,\lambda)=\frac{\sin\rho(a+x-1)}\rho
+O\Big(\frac1{\rho^2}\exp(|{\rm Im}\rho|(|a|+1-x))\Big),
 \quad \lambda\to\infty,
\end{equation}
uniformly with respect to $x\in[0,1],$ which along with the classical
 asymptotics
$$
S(x,\lambda)=\frac{\sin\rho x}\rho +O\Big(\frac1{\rho^2} \exp(|{\rm Im}|x)\Big), \quad \lambda\to\infty
$$
and (\ref{a.1.1}) give the asymptotic formula
\begin{equation}\label{a.3}
G(x,t,\lambda)=\frac1{2\lambda\Delta(\lambda)}\left\{
\begin{array}{l} \cos\rho(a+x-1+t) -\cos\rho(a+x-1-t)\\[3mm]
\qquad\qquad\displaystyle +O\Big(\frac1\rho\exp(|{\rm Im}\rho|(|a|+1-x+t))\Big), \quad t\le x,\\[3mm]
\cos\rho(a+x-1+t) -\cos\rho(a+t-1-x)\\[3mm]
\qquad\qquad\displaystyle +O\Big(\frac1\rho\exp(|{\rm Im}\rho|(|a|+1-t+x))\Big), \quad t\ge x.
\end{array}\right.
\end{equation}

Consider the set
\begin{equation}\label{a.3.1}
D_\varepsilon(\sigma):=\Big\{\lambda=\rho^2: \Big|\rho-\frac{\pi k}\sigma\Big|\ge\varepsilon, k\in{\mathbb
Z}\Big\}, \quad \varepsilon>0, \quad \sigma\ne0.
\end{equation}

\medskip
 {\bf Example A1.} Let $a\le0.$ Then the problem $R(a,q)$ is Birkhoff
 regular for any $q(x).$ Indeed, according to (\ref{3a}), we have the
 estimate
 $$
 |\Delta(\lambda)|\ge\frac{C_\varepsilon}{|\rho|}\exp(|{\rm Im}\rho|(1-a)),
 \quad \lambda\in D_\varepsilon(1-a), \quad |\lambda|\ge r_\varepsilon,
 $$
 for sufficiently large $r_\varepsilon,$ which along with (\ref{a.3})
 give estimate (\ref{a.1}) for $\theta=-1/2.$

\medskip
{\bf Example A2.} For any $a>0,$ the problem $R(a,0)$ is not regular even in the Stone sense. Indeed, for the
zero potential, the Green's function has the form
$$
G(x,t,\lambda)=\frac1{\rho\sin\rho(a-1)}\left\{
\begin{array}{l} \sin\rho(a+x-1)\sin\rho t, \quad t\le x,\\[3mm]
\sin\rho(a+t-1)\sin\rho x, \quad t\ge x,
\end{array}\right.
$$
as soon as $a\ne1,$ while it does not exist for $a=1,$ because in this case $\Delta(\lambda)\equiv0.$ Hence,
the problem $R(1,0)$ is automatically irregular. For $a\ne1,$ the latter representation implies the following
estimates:
$$
|G(1,1,\lambda)|\ge\frac{C_\varepsilon}{|\rho|}\exp(2|{\rm Im}\rho|a), \quad \lambda\in D_\varepsilon(1)\cap
D_\varepsilon(a), \quad a\in(0,1)\cup[3,\infty),
$$
$$
\Big|G\Big(\frac{a-1}2,\frac{a-1}2,\lambda\Big)\Big|\ge\frac{C_\varepsilon}{|\rho|}\exp(|{\rm Im}\rho|(a-1)),
\quad \lambda\in D_\varepsilon\Big(\frac{a-1}2\Big)\cap D_\varepsilon\Big(\frac{3a-3}2\Big), \quad a\in(1,3),
$$
which imply impossibility of estimate (\ref{a.1}) for any finite $\theta.$

\medskip
{\bf Example A3.} Any problem $R(1,q)\in{\mathfrak R}$ is Stone regular. Indeed, since (\ref{a.3}) takes the
form
$$
G(x,t,\lambda)=\frac1{2\lambda\Delta(\lambda)}\Big(\cos\rho(x+t) -\cos\rho(x-t) +O\Big(\frac1\rho\exp(|{\rm
Im}\rho|(2-|x-t|))\Big)\Big)
$$
and formula (\ref{partsDelta}) implies the estimate
$$
|\Delta(\lambda)|\ge\frac{C_\varepsilon}{|\rho|^3}\exp(2|{\rm Im}\rho|), \quad \lambda\in D_\varepsilon(2),
\quad |\lambda|\ge r_\varepsilon,
$$
we arrive at (\ref{a.1}) with $\theta=1/2.$
\\

\noindent{\large\bf Appendix B}
\\

In what follows, we let $q(x)$ be an arbitrary real-valued function in $L_2(0,1),$ and consider the following
inverse problem.

\medskip
{\bf Inverse Problem B1.} Given the spectrum $\{\lambda_k\}$ of the problem $R(a,q),$ find $q(x).$

\medskip
We show that the solution of Inverse Problem~B1 is unique if and only if $a\in(-\infty,-1]\cup(1,\infty)$ and
obtain necessary and sufficient conditions of its solvability for all real $a.$ The case $a=1$ is exceptional
and treated separately in the following theorem.

\medskip
{\bf Theorem B1. }{\it For any sequence $\{\lambda_k\}$ of complex numbers to be the spectrum of the boundary
value problem $R(1,q)$ with a real-valued square-integrable potential $q(x),$ it is necessary and sufficient
that the infinite product in (\ref{3}) is convergent for each complex $\lambda$ and the corresponding
function $\Theta(\lambda)$ has the form
\begin{equation}\label{b.1}
\Theta(\lambda) =\int_0^2 u(x) \frac{\sin\rho x}\rho\, dx, \quad u(x)\in W_{2,{\mathbb R}}^1[0,1], \quad
u(2)=0,
\end{equation}
where $W_{2,{\mathbb R}}^1[0,1]$ is the real version of the space $W_2^1[0,1].$}

\medskip
{\it Proof.} For the necessity, it is sufficient to note that, according to the first equality in (\ref{3})
along with (\ref{Delta}), we have
\begin{equation} \label{b.2}
\Theta(\lambda) = \frac\alpha{\rho^2} +\int_0^2 g(t) \frac{\cos \rho t}{\rho^2} \, dt,
\end{equation}
where
$$
\alpha = -\frac\omega{2\gamma}, \quad g(t) =\frac{w(t)}\gamma, \quad \int_0^2 w(t)\, dt =\frac\omega2,
$$
Thus, integrating by parts in (\ref{b.2}), we get (\ref{b.1}) with
$$
u(x)=-\int_x^2 g(t)\, dt.
$$

Let us prove the sufficiency. Integrating by parts in (\ref{b.1}) we obtain (\ref{b.2}) with $g(x)=u'(x)$ and
$\alpha=u(0).$ Calculate the functions $g_0(x)$ and $g_1(x)$ by the formula
$$
g_j(x)=g(1+x)+(-1)^jg(1-x), \quad x\in(0,1), \quad j=0,1,
$$
For $j=0,1$ and  $\gamma\in{\mathbb R}\setminus\{0\},$ denote by $\{\lambda_{k,j}\}_{k\ge1}$ the sequence of
zeros (with an account of multiplicity) of the function $\Delta_j(\lambda),$ determined by the corresponding
formula in (\ref{1-3}) with $\omega=-\gamma\alpha$ and $w_j(t)=\gamma g_j(t).$ Hence, the asymptotics
(\ref{asympt}) holds and, according to the proof of Theorem~1 in \cite{BYY}, the zeros interlace as in
(\ref{interlace}) for sufficiently small positive $|\gamma|.$ Then, by virtue of Theorem~2, there exists a
real-valued potential $q_\gamma(x)\in L_2(0,\pi)$ such that $\{\lambda_{k,j}\}_{k\ge1}$ is the spectrum of
the problem ${\cal L}_j(q_\gamma)$ for $j=0,1.$ Moreover, as in Section~2, one can show that the
characteristic function of the problem $R(1,q_\gamma)$ coincides with $\gamma\Theta(\lambda).$ $\hfill\Box$

\medskip
We note that, as was first established in \cite{AGP}, the constructed potential $q_\gamma(x)$ is uniquely
determined by fixing the value $\gamma.$ Otherwise, there are infinitely many potentials corresponding to one
and the same spectrum $\{\lambda_k\}$ (see \cite{BYY}), although not any $\gamma$ may lead to some potential
$q_\gamma(x)$ (see \cite{YB20}). However, specification of the spectrum of $R(a,q)$ does determine the
potential uniquely as soon as $a>1$ or $a\le-1,$ which can be seen in the proof of the next theorem.

\medskip
{\bf Theorem B2. }{\it For any sequence $\{\lambda_k\}$ of complex numbers to be the spectrum of the boundary
value problem $R(a,q)$ with $a\in(-\infty,-1]\cup(1,\infty)$ and a real-valued square-integrable potential
$q(x),$ it is necessary and sufficient that the following two conditions are fulfilled:

(i) The infinite product in (\ref{3}) is convergent for each complex $\lambda$ and the corresponding function
$\Delta(\lambda)$ with
\begin{equation}\label{b.3}
\gamma=2\frac{1-a}\pi
\lim_{n\to\infty}\frac{(-1)^{n+1}}{2n-1}\Big(\Theta\Big(\frac{\pi^2}{(1-a)^2}\Big(n-\frac12\Big)^2\Big)\Big)^{-1}
\end{equation}
has the form (\ref{3a}) with some real-valued function $w(t)\in L_2(a-1,a+1)$ and $\omega\in{\mathbb R};$

(ii) Zeros of the functions $\Delta_0(\lambda)$ and $\Delta_1(\lambda)$ constructed by (\ref{1-3}) with
$w_0(t)$ and $w_1(t)$ determined by (\ref{wa}), i.e.
\begin{equation}\label{b.3.1}
w_j(t)=w(a+t)+(-1)^jw(a-t), \quad t\in(0,1), \quad j=0,1,
\end{equation}
are real and interlacing.}

\medskip
{\it Proof.} By necessity, both the representations (\ref{3}) and (\ref{3a}) are already established. Thus,
for~(i), it is sufficient to prove (\ref{b.3}). According to (\ref{3a}), for any fixed positive $\varepsilon$
we have
\begin{equation}\label{b.4}
\frac{\rho\Delta(\lambda)}{\sin\rho(1-a)}\to1, \quad \lambda\to\infty, \quad \lambda\in D_\varepsilon(1-a),
\end{equation}
where the set $D_\varepsilon(\sigma)$ is determined in (\ref{a.3.1}). In particular, we have
$\{\eta_n\}_{n\in{\mathbb N}}\subset D_\varepsilon(1-a)$ as soon as $\varepsilon\in(0,\pi(a-1)/2),$ where
\begin{equation}\label{b.5}
\eta_n=\theta_n^2, \quad \theta_n=\frac\pi{1-a}\Big(n-\frac12\Big).
\end{equation}
Thus, (\ref{b.4}) implies
$$
(-1)^{n+1}\theta_n\Delta(\eta_n)\to1, \quad n\to\infty,
$$
which along with (\ref{b.5}) and the first identity in (\ref{3}) imply (\ref{b.3}). For~(ii), it remains to
note that the functions $\Delta_0(\lambda)$ and $\Delta_1(\lambda)$ are the characteristic functions of the
problems ${\cal L}_0(q)$ and~${\cal L}_1(q),$ respectively. According to Theorem~2, their zeros are real,
simple and interlace.

Let us prove the sufficiency of the conditions (i) and (ii). For $j=0,1,$ we let $\{\lambda_{k,j}\}_{k\ge1}$
be zeros of the function $\Delta_j(\lambda)$ determined in (ii). Then they have asymptotics (\ref{asympt})
and, hence, the interlacement implies (\ref{interlace}). According to Theorem~2, there exists a real-valued
potential $q(x)\in L_2(0,1)$ such that the functions $\Delta_0(\lambda)$ and $\Delta_1(\lambda)$ are the
characteristic functions of the problems ${\cal L}_0(q)$ and ${\cal L}_1(q),$ respectively. Consider the
problem $R(a,q)$ and let $\tilde\Delta(\lambda)$ be its characteristic function, which possesses the
representation:
$$
\tilde\Delta(\lambda)=\frac{\sin\rho(1-a)}{\rho}-\tilde\omega\frac{\cos\rho(1-a)}{2\rho^2}+ \int_{a-1}^{a+1}
\tilde w(t)\frac{\cos\rho t}{\rho^2}\,dt,
$$
where
$$
\tilde\omega=\int\limits_0^1q(x)\,dx, \quad \tilde w(t) = \frac{1}{2} \left\{\begin{array}{cc}
w_0(a-t)-w_1(a-t), & t\in(a-1,a),\\[2mm]
w_0(t-a)+w_1(t-a), & t\in(a,a+1).
\end{array} \right.
$$
Comparing these formulae with (\ref{omega}) and (\ref{wa}), respectively, according to the representation
(\ref{3a}), we arrive at $\tilde\Delta(\lambda)\equiv\Delta(\lambda)$ and hence the sequence $\{\lambda_k\}$
is the spectrum of $R(a,q).$ $\hfill\Box$

\medskip
Finally, we consider the case $a\in(-1,1),$ when Inverse Problem~B1 is not uniquely solvable. Indeed,
according to representation (\ref{3a}), specification of $\Delta(\lambda)$ determines only the even part
$w_+(t)$ of the function $w(t)$ on the interval $(-b,b),$ where
$$
b=\min\{1-a,1+a\}, \quad w_+(t)=\frac{w(t)+w(-t)}2.
$$
Hence, (\ref{3a}) takes the form
\begin{equation}\label{b.5.1}
\Delta(\lambda)=\frac{\sin\rho(1-a)}{\rho}-\omega\frac{\cos\rho(1-a)}{2\rho^2}+  2\int_0^b
w_+(t)\frac{\cos\rho t}{\rho^2}\,dt +\int_{a_1}^{a_2} w(t)\frac{\cos\rho t}{\rho^2}\,dt,
\end{equation}
where $a_1={\rm sgn}(a)(1-a)$ and $a_2={\rm sgn}(a)(1+a).$ Therefore, in order to recover $w(t)$ completely
and, thus, to fix a unique solution of Inverse Problem~B1, one should additionally specify the odd part
$w_-(t)$ of $w(t)$ on $(-b,b):$
$$
w_-(t)=\frac{w(t)-w(-t)}2.
$$
Analogously to Theorem~B2, one can prove the following theorem, which gives necessary and sufficient
conditions for solvability (not unique) of Inverse Problem~B1 when $a\in(-1,1).$

\medskip
{\bf Theorem B3. }{\it For any sequence $\{\lambda_k\}$ of complex numbers to be the spectrum of the problem
$R(a,q)$ with $a\in(-1,1)$ and a real-valued square-integrable potential $q(x),$ it is necessary and
sufficient that, besides condition (i) in Theorem~B2, the following condition is fulfilled:

(ii') There exists a real-valued function $w_-(t)\in L_2(0,b)$ such that zeros of the functions
$\Delta_0(\lambda)$ and $\Delta_1(\lambda),$ constructed by (\ref{1-3}) with $w_0(t)$ and $w_1(t)$ determined
by (\ref{b.3.1}), are real and interlacing. Here $w(t)$ is determined on $(a_1,a_2)$ by the relation
(\ref{b.5.1}), while on $(-b,b)$ it is determined by the formula
$$
w(t)=\left\{\begin{array}{l}
w_+(t)-w_-(t), \quad t\in(-b,0),\\[3mm]
w_+(t)+w_-(t), \quad t\in(0,b).
\end{array}\right.
$$}

{\bf Remark B1.} For $a=0,$ it is well-known that Inverse Problem~B1 is solvable if and only if the numbers
$\lambda_k$ are real, simple and obey the asymptotics
$$
\lambda_k= \pi^2k^2+\omega+\varkappa_k, \quad \{\varkappa_k\}\in l_2, \quad k\ge1.
$$
Meanwhile, under such assumptions, conditions (i) and (ii') can be checked separately. Indeed, the fulfilment
of (i) can be proved by using Lemma~3.3 in \cite{But07}, while (ii') follows from Theorem~2.

\medskip
{\bf Remark B2.} Conditions (ii) and (ii') in Theorems~B2 and~B3 can be formulated in terms of a Nevanlinna
function. By definition, a complex function belongs to the Nevanlinna class, if it is analytic on the open
upper half-plane and has non-negative imaginary part there. Consider the meromorphic function
$$
M(\lambda):= \frac{\Delta_0(\lambda)}{\Delta_1(\lambda)},
$$
where the functions $\Delta_0(\lambda)$ and $\Delta_1(\lambda)$ are determined by formulae (\ref{1-3}) with
some real number~$\omega$ and real-valued square-integrable functions $w_0(x)$ and $w_1(x).$ Thus, each of
conditions (ii) and (ii') is equivalent to belonging of the function $M(\lambda)$ to the Nevanlinna class.
Indeed, for the latter it is necessary and sufficient that zeros and poles of $M(\lambda)$ interlace, which
can be proved analogously to Theorem~1 on page~308 in \cite{Lev}. For convenience of the reader, we provide
the crucial arguments of the proof. First of all, we note that the following representations hold:
$$
\Delta_j(\lambda)=\prod_{k=1}^\infty\frac{\lambda_{k,j}-\lambda}{\pi^2(k-j/2)^2}, \quad j=0,1,
$$
(see, e.g., \cite{FY}). Assume that (\ref{interlace}) holds and let $N\in{\mathbb N}$ be such that
$\lambda_{N-1,0}<0<\lambda_{N+1,1},$ where $\lambda_{0,0}=-\infty.$ Then we have the formulae
$$
M(\lambda)=C\frac{\lambda-\lambda_{N,0}}{\lambda-\lambda_{N,1}} \prod_{k\ne
N}\frac{1-\frac\lambda{\lambda_{k,0}}}{1-\frac\lambda{\lambda_{k,1}}}, \quad C=\Big(1-\frac1{2N}\Big)^2
\prod_{k\ne N}\frac{\lambda_{k,0}}{\lambda_{k,1}}\Big(1-\frac1{2k}\Big)^2.
$$
Since the formula for $C$ possesses an even number of negative multipliers, we have $C>0.$ Moreover, since
$\lambda_{k,0}\lambda_{k,1}>0$ as soon as $k\ne N,$ we obtain
$$
\arg M(\lambda)=\sum_{k=1}^\infty d_n \in(0,\pi), \quad d_n:=\arg(\lambda-\lambda_{k,0})
-\arg(\lambda-\lambda_{k,1})>0, \quad {\rm Im}\lambda>0,
$$
which proves the sufficiency. Let $M(\lambda)$ now be a Nevanlinna function, i.e. $\arg M(\lambda)\in[0,\pi]$
for ${\rm Im}\lambda>0$ and, symmetrically,  $\arg M(\lambda)\in\{\pi\}\cup(-\pi,0]$ for ${\rm Im}\lambda<0,$
since $\overline{M(\overline\lambda)}=M(\lambda).$ Then all its zeros and poles should be real. Otherwise, a
circuit around any single zero or pole lying in the open upper or lower half-plane would increment $\arg
M(\lambda)$ by not less than $2\pi,$ which is impossible. By the same means, we establish that all zeros and
poles of $M(\lambda)$ are simple, and on any interval their numbers differ by no more than one, i.e. zeros
and poles interlace.

Finally, it can be reminded that $M(\lambda)$ is the Weyl function of the operator generated by the
differential expression $\ell$ and the boundary conditions $y(0)=y'(1)=0$ (see, e.g., \cite{FY}).
\\

{\bf Funding.} The first and the third authors were supported by
Russian Foundation for Basic Research (Project No. 20-31-70005).
 The second author is  supported
 by CONACYT Project A1-S-31524 and CIC-UMSNH, Mexico.

\end{document}